\documentclass{amsart}
\usepackage{amssymb}

\theoremstyle{plain}
\newtheorem{thm}{Theorem}[section]
\newtheorem{lem}[thm]{Lemma}
\newtheorem{prop}[thm]{Proposition}

\theoremstyle{definition}

\newcommand{\Cplx}{{\bf{C}}}

\newcommand{\Zed}{{\bf{Z}}}
\newcommand{\GL}{{\mathrm{GL}}}

\newcommand{\suchthat}{\thinspace | \thinspace}

\newcommand{\co}{\colon\thinspace}

\title{Generalized Long-Moody representations of braid groups}

\author{Stephen Bigelow}
\address{Department of Mathematics, University of California at Santa Barbara,
         CA 93106}
\email{bigelow@math.ucsb.edu}

\author{Jianjun Paul Tian}
\address{Mathematics Department, the College of William and Mary, Williamsburg, VA 23187}
\email{jptian@math.wm.edu}

\thanks{Dedicated to the memory of Professor Xiao-Song Lin}

\begin{document}

\begin{abstract}
Long and Moody give a method of constructing representations of the
braid groups $B_n$. We discuss some ways to generalize their
construction. One of these gives representations of subgroups of
$B_n$, including the Gassner representation of the pure braid group
as a special case. Another gives representations of the Hecke
algebra.
\end{abstract}

\maketitle

\section{Introduction}
\label{sec:intro}

The braid group $B_n$ is the group with presentation
$$\langle \sigma_1,\dots,\sigma_{n-1} \suchthat
  \sigma_i \sigma_j = \sigma_j \sigma_i \mbox{ if } |i-j| > 1, \thickspace
  \sigma_i \sigma_j \sigma_i = \sigma_j \sigma_i \sigma_j \mbox{ if } |i-j|=1
  \rangle.
$$
It has a rich and interesting representation theory
\cite{birman} \cite{formanek}.

D. D. Long, in joint work with J. Moody, gave a construction that
inputs a representation of a certain semidirect product $F_n \rtimes
B_n$ and outputs a new representation of $B_n$. At first this
construction may not sound useful, since its input seems to be
harder to come by than its output. However there are at least two
factors in its favor. First, the output seems
to be in some sense more sophisticated than the input. For example,
one can start with a one dimensional representation and obtain the
Burau representation. Second, the semidirect product $F_n \rtimes
B_n$ is a subgroup of $B_{n+1}$. Thus one can recursively apply the
algorithm to construct a sequence of increasingly complicated
representations of $B_n$ for all $n$.

Let $F_n$ be the free group $\langle g_1,\dots,g_n \rangle$. There
is a well known action of $B_n$ on $F_n$, which gives rise to a
semidirect product $F_n \rtimes B_n$. This requires an arbitrary
choice of convention. In this paper, the relations are as
follows.
\begin{itemize}
\item $g_{i+1} \sigma_i = \sigma_i g_i$,
\item $g_i \sigma_i = \sigma_i g_i g_{i+1} g_i^{-1}$,
\item $g_j \sigma_i = \sigma_i g_j$ if $j \not \in \{i,i+1\}$.
\end{itemize}

The following is Theorem 2.1 in \cite{dL94}, which Long credits to
be a joint work with Moody.

\begin{thm}
\label{thm:lm}
Given a representation $\rho \co F_n \rtimes B_n \to \GL(V)$
we may construct a representation
$\rho^+ \co B_n \to \GL(V^{\oplus n})$.
\end{thm}

The method of construction can be summarized as follows. The
representation $\rho$ determines a system of local
coefficients on the $n$-times punctured disk $D_n$. The first
cohomology of $D_n$ with these local coefficients then turns out to
be isomorphic to $V^{\oplus n}$. The braid group $B_n$ is the
mapping class group of $D_n$, so it has a natural induced action on
this vector space. For details, see \cite{dL94}.

Long also describes how to explicitly compute the
matrices for $\rho^+$ in terms of the matrices for
$\rho$. In this paper we will only use this more concrete (but less
well motivated) definition.

In Section \ref{sec:lm}, we give a proof of Theorem \ref{thm:lm}. In
Section \ref{sec:subgroup}, we generalize the Long-Moody
construction to subgroups of $B_n$. We then discuss important
special cases of this construction. In Section \ref{sec:gassner}, we
apply the construction to the pure braid group $P_n$, obtaining the
Gassner representation as a special case. In Section
\ref{sec:lawrence}, we describe how these methods can be used to
obtain Lawrence's representations. These are of interest because
they contain all representations of the Hecke algebra. Finally, in
Section \ref{sec:reduced}, we give the ``reduced'' version of the
Long-Moody construction.

For convenience, we will work with finite dimensional
representations over the field $\Cplx$. The construction works
equally well for infinite dimensional representations, and for
representations over any field or commutative ring.

Note that, while the main results of this paper are stated as
Theorems and Corollaries and so on, they are really {\em
constructions}. Taken literally, a statement of the form ``Given one
representation, we may construct another representation'' is
trivial. Of course, the interesting thing is the method of
construction.

\section{The universal Long-Moody representation}
\label{sec:lm}

The aim of this section is to give a proof of Theorem \ref{thm:lm}.
Our proof is not exactly new, but it does make explicit an idea that
is implicit in \cite{dL94}. This is the idea of a representation we
call the ``universal Long-Moody representation''.

\begin{proof}[Proof of Theorem \ref{thm:lm}]
Recall that we are given a representation $\rho \co F_n \rtimes B_n
\to \GL(V)$, and must construct a representation $\rho^+ \co B_n \to
\GL(V^{\oplus n})$. First we define the ``universal Long-Moody
representation'' of $B_n$. This is a homomorphism from $B_n$ to the
group of invertible $n$ by $n$matrices
with entries in the group ring of $F_n \rtimes B_n$.

Let $F_n \rtimes B_n$ be the semidirect product as defined earlier.
Let $\Zed[F_n \rtimes B_n]$ denote the group ring of this group. For
$i=1,\dots,n-1$, let $R_i$ be the following two by two matrix with
entries in $\Zed[F_n \rtimes B_n]$.
$$
R_i =
\left[
\begin{array}{cc}
0 &g_i    \\
1 &1-g_i
\end{array}
\right].
$$
This is invertible.
$$
R_i^{-1} =
\left[
\begin{array}{cc}
(1-g_i^{-1}) & 1\\
g_i^{-1}     & 0
\end{array}
\right].
$$

Let $\phi$ be the following homomorphism from $B_n$ to the group of
invertible $n$ by $n$ matrices with entries in $\Zed[F_n \rtimes
B_n]$.
$$
\phi(\sigma_i) =
\sigma_i
\left[
\begin{array}{ccc}
I_{i-1}               \\
       &R_i           \\
       &   &I_{n-i-1}
\end{array}
\right].
$$
Here, $I_{i-1}$ and $I_{n-i-1}$ are the appropriately sized identity
matrices, and $\sigma_i$
acts on the matrix by scalar multiplication on the left.

It remains to check that $\phi$ satisfies the braid relations.
The most difficult is the braid relation
$$\phi(\sigma_i    ) \phi(\sigma_{i+1}) \phi(\sigma_i    ) =
  \phi(\sigma_{i+1}) \phi(\sigma_i    ) \phi(\sigma_{i+1})$$
One can verify this by direct calculation, using the relations in
$F_n \rtimes B_n$.

Having defined $\phi$, it is now easy to construct $\rho^+$. Choose
a finite basis for $V$. (This is not strictly necessary - we leave
it to the reader to fill in the details when $V$ is infinite
dimensional.) The image of $\rho$ now consists of square matrices.
For any $\beta \in B_n$, let $\rho^+(\beta)$ be the block matrix
obtained by applying $\rho$ to each entry of $\phi(\beta)$. This
is an element of $\GL(V^{\oplus n})$, as required.
\end{proof}

The simplest nontrivial example of this construction is
stated in the following proposition.

\begin{prop}
\label{prop:burau} If $\rho \co F_n \rtimes B_n \to \Cplx^*$ is a
one dimensional representation then $\rho^+ \co B_n \to
\GL(n,\Cplx)$ is the unreduced Burau representation, up to
specializing and rescaling.
\end{prop}

\begin{proof}
Suppose $\rho \co F_n \rtimes B_n \to \Cplx^*$ is a one dimensional
representation. The generators $\sigma_i$ of $B_n$ are all conjugate
to each other, so must all be mapped to the same value $s \in
\Cplx^*$. Similarly, the generators $g_i$ are all conjugate to each
other in $F_n \rtimes B_n$, so must all be mapped to the same value
$t \in \Cplx^*$. Then
$$
\rho^+(\sigma_i) =
s\left[
\begin{array}{cccc}
I_{i-1}               \\
       &0&t           \\
       &1&1-t         \\
       & &   &I_{n-i-1}
\end{array}
\right].
$$
This is the familiar matrix of the unreduced Burau representation,
rescaled by a factor of $s$.
\end{proof}

\section{The Long-Moody construction for subgroups of $B_n$.}
\label{sec:subgroup}

Suppose $G$ is a subgroup of $B_n$. Then $F_n \rtimes G$ is a
subgroup of $F_n \rtimes B_n$. The aim of this section is to prove
the following theorem.

\begin{thm}
\label{thm:subgroup} Given a subgroup $G$ of $B_n$ and
a representation $\rho \co F_n \rtimes G \to \GL(V)$,
we may construct a representation $\rho^+ \co G \to \GL(V^{\oplus n})$.
\end{thm}

The construction is the same as that of Long and Moody. Given an
element $\beta \in G$, we let $\rho^+(\beta)$ be the result of
applying $\rho$ to each of the entries in the matrix $\phi(\beta)$.
We need only check that these entries lie in the domain of $\rho$.
Thus it remains only to
prove the following lemma.

\begin{lem}
\label{lem:subgroup} If $\beta \in B_n$ then every entry of the
matrix $\phi(\beta)$ is of the form $a \beta$ for some $a \in
\Zed[F_n]$.
\end{lem}

\begin{proof}
Write $\beta$ as a word in the generators $\sigma_i$.
Recall that $\phi(\sigma_i)$ is
the scalar $\sigma_i$ times a matrix with entries in $\Zed[F_n]$.
Thus $\phi(\beta)$ is a product of terms,
each of which is either a scalar $\sigma_i^{\pm 1}$
or a matrix with entries in $\Zed[F_n]$.
The product of the terms $\sigma_i^{\pm 1}$ in this product,
taken in order,
is equal to $\beta$.

Observe that if $A$ is any matrix with entries in $\Zed[F_n]$
and $\sigma_i$ is a generator of $B_n$
then $\sigma_i A$ is equal to $A' \sigma_i$
for some matrix $A'$ with entries in $\Zed[F_n]$.
To see this,
apply the relations of $F_n \rtimes B_n$ to each term in each entry of $A$.
By using this observation repeatedly,
we can move the terms $\sigma_i^{\pm 1}$
to the right of the product expression for $\phi(\beta)$.
We obtain an expression of the form $A \beta$,
where $A$ is a matrix with entries in $\Zed[F_n]$.
\end{proof}

\section{The Gassner representation}
\label{sec:gassner}

Let $P_n$ denote the pure braid group. The aim of this section is to
prove the following theorem.

\begin{thm}
\label{thm:gassner} Given a representation $\rho \co P_{n+1} \to
\GL(V)$ we may construct an $n$ parameter family of representations
$\rho^+_{t_1,\dots,t_n} \co P_n \to \GL(V^{\oplus n})$.
\end{thm}

Long achieved a similar result by different methods in
\cite{dL89} and \cite{dL89-2}. The construction there uses an action of
the braid group by diffeomorphisms on a certain representation
variety. To obtain a linear representation, he takes the derivative
of this action at any point that is fixed by the action of $B_n$, or $P_n$.
Long was most interested in the Gassner representation,
but the approach is quite general, and could give a result similar
to Theorem~\ref{thm:gassner}.

Our construction is more along the lines of \cite{dL94}, especially
the following theorems, which are \cite[Theorem 2.4]{dL94} and
\cite[Corollary 2.6]{dL94} respectively.

\begin{thm}
Given a representation $\rho$ of $B_{n+1}$, we may construct a new
representation $\rho^+$ of $B_n$.
\end{thm}

\begin{thm}
\label{thm:parameter} Given a representation $\rho$ of $F_n \rtimes
B_n$, we may construct a one parameter family of representations
$\rho_t^+$ of $B_n$.
\end{thm}

The first follows immediately from the fact that $F_n \rtimes B_n$
is a subgroup of $B_{n+1}$. The second uses $\rho$ to define a one
parameter family of representations $\rho_t$ of $F_n \rtimes B_n$,
and then applies the Long-Moody construction to these. The
representations $\rho_t$ are just ``rescalings'' of $\rho$, but the
resulting family of representations $\rho^+_t$ is often more
interesting than $\rho^+$.

To prove Theorem \ref{thm:gassner}, we apply the same ideas in the
context of the pure braid group.

\begin{proof}[Proof of Theorem \ref{thm:gassner}]
Let $F_n$ be the subgroup of $P_{n+1}$ consisting of all braids in
which the first $n$ strands are straight. This is the free group
with generators
$$g_i
  = (\sigma_n\dots\sigma_{i+1}) \sigma_i^2
    (\sigma_n\dots\sigma_{i+1})^{-1}.$$
Let $P_n$ be a subgroup of $P_{n+1}$ in the usual way. It is well
known that $P_n$ and $F_n$ together generate $P_{n+1}$. Indeed,
$P_{n+1}$ is isomorphic to the subgroup $F_n \rtimes P_n$, of our
semidirect product $F_n \rtimes B_n$.

Fix $t_1,\dots,t_n \in \Cplx^*$. Let $\rho_{t_1,\dots,t_n} \co
P_{n+1} \to \GL(V)$ be the representation such that
$$\rho_{t_1,\dots,t_n} (g_i)  = t_i \rho(g_i)$$
for any of the generators $g_i$ of $F_n$, and
$$\rho_{t_1,\dots,t_n} (\beta) = \rho(\beta)$$
for any $\beta \in P_n$.

To see that $\rho_{t_1,\dots,t_n}$ is well defined, note that if
$\beta \in P_n$ and $i \in \{1,\dots,n\}$ then
$$\beta g_i \beta^{-1} = w g_i w^{-1}$$
for some $w \in F_n$. These relations determine the semidirect
product structure of $P_{n+1}$, and they are preserved by
$\rho_{t_1,\dots,t_n}$.

Now simply apply Theorem \ref{thm:subgroup} to $\rho_{t_1,\dots,t_n}$
to obtain the desired representation of $P_n$.
\end{proof}

As an example, we compare the results of Theorems
\ref{thm:parameter} and \ref{thm:gassner} when $\rho$ is the trivial
representation. If $\rho$ is the trivial representation of $F_n
\rtimes B_n$ then Theorem \ref{thm:parameter} gives the unreduced
Burau representation $\rho^+_t$ of $B_n$ (by Proposition
\ref{prop:burau}). If $\rho$ is the trivial representation of
$P_{n+1}$ then $\rho^+_{t_1,\dots,t_n}$ is the unreduced Gassner
representation of $P_n$.

For more sophisticated examples, we could think of Theorem
\ref{thm:gassner} as giving a ``Gassner'' version of representations
of $B_n$ other than the Burau representation.

\section{Lawrence's construction.}
\label{sec:lawrence}

The Hecke algebra is a certain quotient of the group algebra $\Cplx
B_n$. For our purposes, a representation of the Hecke algebra is any
representation $\rho$ of $B_n$ such that the matrices
$\rho(\sigma_i)$ all satisfy a quadratic relation. Such
representations are the subject of a large body of ongoing
research \cite{mathas} \cite{stephen}.

The following proposition suggests a connection between the
representation $\phi$ from Section \ref{sec:lm} and representations
of the Hecke algebra.

\begin{prop}
$(\phi(\sigma_i) + \sigma_i g_i I)(\phi(\sigma_i) - \sigma_i I) = 0$.
\end{prop}

This can be checked by direct calculation. However it does not
qualify as a representation of the Hecke algebra, since the
``scalars'' of this quadratic relation lie in a non-commutative
ring.

By rescaling,
it suffices to consider representations of $B_n$ that satisfy
$$(\rho(\sigma_i) + qI)(\rho(\sigma_i) - I) = 0$$
for some $q \in \Cplx^*$.
For generic values of $q$,
such representations of $B_n$
can be enumerated by the partitions of $n$. In
\cite{rL93}, Lawrence gives a topological construction of the
representations of the Hecke algebra that correspond to partitions
of $n$ of the form $(n-m,m)$. This was generalized to all partitions
in \cite{rL96}.

In \cite[Corollary 2.10]{dL94}, Long describes how his construction
can be used to parallel that of \cite{rL93}. The methods in this
section are intended to parallel those of \cite{rL96} (which had not
appeared at the time of \cite{dL94}). We will not prove any concrete
connection. To do that would require a thorough treatment of the
homological details behind the Long-Moody construction. Instead,
this section could be viewed as a self contained description of
a way to construct representations of the Hecke algebra.

Let $B_{n,m}$ denote the set of braids in $B_{n+m}$ such that the
first $n$ nodes at the bottom are connected to the first $n$ nodes
at the top, and (hence) the last $m$ nodes at the bottom are
connected to the last $m$ nodes at the top. This is generated by
$$\sigma_1,\dots,\sigma_{n-1},\sigma_n^2,\sigma_{n+1},\dots,\sigma_{n+m-1}.$$

\begin{lem}
\label{lem:lawrence}
Given a representation $\rho \co B_{n,m+1} \to \GL(V)$,
we may construct a representation
$\rho^+ \co B_{n,m} \to \GL(V^{\oplus (n+m)})$.
\end{lem}

\begin{proof}
Let $B_{n,m}$ be a subgroup of $B_{n,m+1}$ in the obvious way. Let
$F_{n+m}$ be the set of braids in $B_{n,m+1}$ such that the first
$n+m$ strands are straight. This is the free group with generators
$$g_i
  = (\sigma_{n+m}\dots\sigma_{i+1}) \sigma_i^2
    (\sigma_{n+m}\dots\sigma_{i+1})^{-1}.$$
Then $B_{n,m}$ and $F_{n+m}$ together generate a subgroup $F_{n+m}
\rtimes B_{n,m}$ of $B_{n,m+1}$. The result now follows from Theorem
\ref{thm:subgroup}.
\end{proof}

\begin{thm}
\label{thm:lawrence} Given a representation $\rho \co B_m \to
\GL(V)$, we may construct a one parameter family of representations
$\rho^+_q \co B_n \to \GL(V^{\oplus k})$, where $q \in \Cplx^*$ and
$k = (n+m-1)(n+m-2)\dots(n+1)(n)$.
\end{thm}

\begin{proof}
Suppose we are given $\rho \co B_m \to \GL(V)$ and $q \in \Cplx^*$.
Let $\rho' \co B_{n,m} \to \GL(V)$ be the representation such that
\begin{itemize}
\item $\rho'(\sigma_i) = 1$ for $i=1,\dots,n-1$,
\item $\rho'(\sigma_n^2) = q$,
\item $\rho'(\sigma_{n+i}) = \rho(\sigma_i)$ for $i=1,\dots,m-1$.
\end{itemize}
Now apply Lemma \ref{lem:lawrence} to $\rho'$ recursively $m$ times
to obtain the desired representation of $B_n$.
\end{proof}

For example, suppose $\rho$ is the one dimensional representation of
$B_2$ given by $\rho(\sigma_1) = t$, where $t \in \Cplx^*$. Then
$\rho^+_q$ is the Lawrence-Krammer representation
\cite{rL93} \cite{krammer}.

By analogy to \cite{rL96}, we expect the representation $\rho^+_q$
to be especially interesting in when
$$(\rho(\sigma_i) - q^{-1}I)(\rho(\sigma_i) - I) = 0.$$
However Lawrence still has
to do some more work to obtain a representation of the Hecke algebra
as a subrepresentation of $\rho^+_q$. Rather
than try to imitate this in our purely algebraic setting, we give an
alternative construction that gives a representation of the Hecke
algebra directly.

\begin{thm}
Given $m \le n$ and a representation $\rho \co B_m \to \GL(V)$, we
can construct a representation $\rho^+ \co B_n \to \GL(V \oplus
\dots \oplus V)$, where there are $\binom{n}{m}$ copies of $V$ in the
direct sum. If $\rho$ satisfies the relation
$(\rho(\sigma_i) + qI)(\rho(\sigma_i) - I) = 0$
then so does $\rho^+$.
\end{thm}

\begin{proof}
We start by constructing a ``universal'' Lawrence representation.
Let $\langle q \rangle$ denote the infinite cyclic group generated
by a formal variable $q$.
Our universal Lawrence representation will be
a homomorphism $\phi$ from $B_n$ to the group of
invertible $\binom{n}{m}$ by $\binom{n}{m}$ matrices
with entries in $\Zed[\langle q \rangle \times B_m]$.
Enumerate the rows and columns of this matrix
by basis vectors
$e_{a_1,\dots,a_m}$ for $1 \le a_1 < \dots < a_m \le n$.

To make the definition more readable, we first fix some notation.
Let ${\mathbf a}$ denote a sequence of the form $a_1,\dots,a_j$,
where $1 \le a_1 < \dots < a_j < i$. Let ${\mathbf b}$ denote a
sequence of the form $b_1,\dots,b_k$, where $i+1 < b_1 < \dots < b_k
\le n$. In each of the following, assume that $j$ and $k$ are such
that the length of the whole sequence is $n$.
\begin{eqnarray*}
\phi(\sigma_i)(e_{{\mathbf a},{\mathbf b}})       &=&
               e_{{\mathbf a},{\mathbf b}}            \\
\phi(\sigma_i)(e_{{\mathbf a},i  ,{\mathbf b}})   &=&
               e_{{\mathbf a},i+1,{\mathbf b}}        \\
\phi(\sigma_i)(e_{{\mathbf a},i+1,{\mathbf b}})   &=&
              qe_{{\mathbf a},i  ,{\mathbf b}}
        +(1-q) e_{{\mathbf a},i+1,{\mathbf b}}        \\
\phi(\sigma_i)(e_{{\mathbf a},i,i+1,{\mathbf b}}) &=&
   \sigma_{j+1}e_{{\mathbf a},i,i+1,{\mathbf b}}.
\end{eqnarray*}
We must check that $\phi$ satisfies the braid relations, so that it
gives a well defined representation of $B_n$. This is done by
direct computation and checking cases.
We will only describe some special
cases that give a flavor of the general proof.

If $m=0$ then $\phi$ is the trivial representation.

If $m=1$ then $\phi$ is the unreduced Burau representation.

If $m=n$ then $\phi$ maps each braid $\beta$
to the one by one matrix $[\beta]$.

Finally, consider the case $m=2$ and $n=3$. The matrices for $\phi$
with respect to the basis $(e_{2,3}, e_{1,3}, e_{1,2})$ are as
follows.
\begin{eqnarray*}
\phi(\sigma_1) &=&
\left[
\begin{array}{ccc}
1-q & 1 & 0        \\
q   & 0 & 0        \\
0   & 0 & \sigma_1
\end{array}
\right],                 \\
\phi(\sigma_2) &=&
\left[
\begin{array}{ccc}
\sigma_1 &  0  & 0        \\
0        & 1-q & 1        \\
0        &  q  & 0
\end{array}
\right].
\end{eqnarray*}
One can check that these matrices satisfy the braid relation
$ABA=BAB$.

To prove that $\phi$ satisfies the braid relations in general, we
can break down the image into blocks, each of which resembles one of
the above cases.

Having defined $\phi$, we can define the specific representation
$\rho^+$ in the same way as we did for the Long-Moody construction.
First choose a basis for $V$, and an ordering for
our basis vectors $e_{a_1,\dots,a_k}$.
The images of $\rho$ and $\phi$ can then be written as square
matrices. Given $\beta \in B_n$, let $\rho^+(\beta)$ be the block
matrix obtained by applying $\rho$ to each entry of $\phi(\beta)$
and substituting a value $q \in \Cplx^*$.

By some further computation and checking of cases,
every nonzero entry of the matrix
$$(\phi(\sigma_i) + qI)(\phi(\sigma_i) - I)$$
is of the form
$$(q+\sigma_j)(1-\sigma_j)$$
for some $j$.
Now suppose $\rho$ satisfies the relation
$$(\rho(\sigma_i) + qI)(\rho(\sigma_i) - I) = 0.$$
If we apply $\rho$ to
each entry of $(\phi(\sigma_i)+qI)(\phi(\sigma_i)-I)$ then we
obtain the zero matrix. Thus $\rho^+$ satisfies the same quadratic
relation as $\rho$.
\end{proof}

\section{The reduced Long-Moody construction.}
\label{sec:reduced}

The unreduced Burau representation obtained in Proposition
\ref{prop:burau} is reducible. It is the direct sum of a one
dimensional representation with the reduced Burau representation. In
\cite[Theorem 2.11]{dL94}, Long shows that this happens in general.
We record here how to compute the ``reduced'' version of the
Long-Moody construction.

For $i=1,\dots,n-1$,
let $S_i$ be the following three by three matrix
with entries in $\Zed[F_n \rtimes B_n]$.
$$
S_i =
\left[
\begin{array}{ccc}
1 &  g_i & 0 \\
0 & -g_i & 0 \\
0 &   1  & 0
\end{array}
\right].
$$
This is invertible.
$$
S_i^{-1} =
\left[
\begin{array}{ccc}
1 &   1       & 0 \\
0 & -g_i^{-1} & 0 \\
0 &  g_i^{-1} & 1
\end{array}
\right].
$$

Let $\phi_r$ be the following homomorphism from $B_n$ to the group
of invertible $n-1$ by $n-1$ matrices with entries in $\Zed[F_n
\rtimes B_n]$.
$$
\phi_r(\sigma_i) =
\sigma_i
\left[
\begin{array}{ccc}
I_{i-2}             \\
       &S           \\
       & &I_{n-i-2}
\end{array}
\right].
$$
If $i$ is $1$ or $n-1$ then part of $S$ is ``cut off'' by the edge
of the matrix. In other words,
$$
\phi_r(\sigma_1) =
\sigma_1
\left[
\begin{array}{ccc}
-g_1&0              \\
  1 &1              \\
    & &I_{n-3}
\end{array}
\right],
$$
$$
\phi_r(\sigma_{n-1}) =
\sigma_{n-1}
\left[
\begin{array}{ccc}
I_{n-3}             \\
        &1&g_{n-1}  \\
        &0&-g_{n-1}
\end{array}
\right].
$$


\newcommand{\etalchar}[1]{$^{#1}$}
\providecommand{\bysame}{\leavevmode\hbox to3em{\hrulefill}\thinspace}
\providecommand{\MRhref}[2]{%
  \href{http://www.ams.org/mathscinet-getitem?mr=#1}{#2}
}
\providecommand{\href}[2]{#2}

\end{document}